\documentclass[12pt]{article}
\usepackage{amsmath,amsfonts,amssymb,rotating,amsthm, bm}
\usepackage{hyperref}
\usepackage{array}
\usepackage{breqn}
\usepackage{color}
\usepackage{fancyhdr}
\usepackage{mathptmx}
\usepackage{enumitem}
\usepackage{dsfont}
\usepackage{caption}
\usepackage{subcaption}
\usepackage{booktabs}
\usepackage{float}
\usepackage {extarrows}
\usepackage{verbatim}
\usepackage{arydshln}

\usepackage{youngtab}
\usepackage{ytableau}
\ytableausetup
{mathmode, boxsize=1.5em}

\usepackage{tikz}
\usetikzlibrary{cd}

\usepackage{times}

\allowdisplaybreaks
\def\qed{\nopagebreak\hfill{\rule{4pt}{7pt}}}
\def\proof{\noindent {\it{Proof.} \hskip 2pt}}
\newtheorem{theorem}{Theorem}[section]

\newtheorem{lemma}[theorem]{Lemma}

\numberwithin{equation}{section}

\begin{document}

\begin{center}
{\Large\bf The
Wide Band Cayley Continuants}

William Y.C. Chen$^1$ and Elena L. Wang$^2$

Center for Applied Mathematics\\
Tianjin University\\
Tianjin 300072, P.R. China

\vskip 3mm

Emails: { $^1$chenyc@tju.edu.cn, $^{2}$ling\_wang2000@tju.edu.cn}

\end{center}

\begin{abstract}
The Cayley continuants are referred to the determinants
of tridiagonal matrices in connection with
the Sylvester continuants.    Munarini-Torri found a striking
combinatorial interpretation
of the Cayley continuants in terms of the joint
distribution of the number of odd cycles and the
number of even cycles of permutations of $[n]=\{1,2,\ldots, n\}$.
In view of a general setting,
$r$-regular cycles (with length not
divisible by $r$) and $r$-singular cycles
(with length divisible by $r$) have been
extensively studied
largely related to roots of permutations. We
introduce
the wide band Cayley continuants as an extension of the
original Cayley continuants,
and we show that they  can be interpreted
in terms of the
joint distribution of the number of $r$-regular cycles and the
number of $r$-singular cycles over permutations of $[n]$.

\end{abstract}

\noindent{\bf Keywords:} $r$-regular cycles, $r$-singular
cycles,   roots of permutations, the Cayley continuants.

\noindent{\bf AMS Classification:} 	05A05 
, 15A15 

\section{Introduction}

As a way to compute
Sylvester's continuants,
Cayley \cite{Cayley-1858} introduced the
tridiagonal determinants, as called the Cayley continuants by
 Munarini and Torri
\cite{MT-2005}.
Cayley's approach may be perceived
as a diagonalization technique in the ternomilogy of today,
that is, generalizing first, then specializing.
Such an understanding has been spelled out by
in \cite{MT-2005}. To be precise, for $n\geq 2$, the Cayley continuants
are referred to the $n \times n$ tridiagonal determinants,
where the zero entries are left blank,
\begin{align*}
U_n(x,y)=\begin{vmatrix}
  x &  1&  &  &  &  \\
  y & x &  2&  & &  \\
   & y-1 &  x &  \ddots &   & \\
  &  & \ddots  &  \ddots &    \ddots &   \\
  & &  & \ddots &  x   &   n-1\\
  &  &  &    & y-n+2&  x
\end{vmatrix}_{n \times n}.
\end{align*}
For $n=0,1$, the initial values are given by
$U_0(x,y)=1$ and $U_1(x,y)=x$.
The determinants satisfy the three-term recurrence for $n \ge 2$,
$$U_{n}(x,y)=x U_{n-1}(x,y)
-(n-1)(y -n+2)U_{n-2}(x,y).$$
From the above recursion, Cayley \cite{Cayley-1858} derived the
exponential generating function of $U_n(x,y)$ and
he further computed the Sylvester's determinants. An exposition
of Cayley's approach can be found in \cite{MT-2005}.
It is noteworthy that several classical
polynomials can be derived from the
Cayley continuants
such as the Meixner
polynomials of the first kind and the Mittag-Leffler polynomials,
see also \cite{MT-2005}.

A striking combinatorial
interpretation of the Cayley continuants
has been found by Munarini and Torri \cite{MT-2005}
in terms of the joint
distribution of the number of
odd cycles and the number of even cycles.
Munarini \cite{Munarini-2022} obtained a representation of
the Cayley continuants by means of the umbral operators.

For $n \geq 1$, let $S_n$ be the
set of permutations of $[n]=\{1,2,\ldots, n\}$.
For a permutation $\sigma \in S_n$,
let $o(\sigma)$ and $e(\sigma)$ be
the number of odd cycles and the
number of even cycles of $\sigma$,
respectively. Then for $n \geq 1$,
\begin{align*}
    U_n(x,y)=\sum_{\sigma \in S_n}
    x^{o(\sigma)} (-y)^{e(\sigma)}.
\end{align*}

To avoid the minus signs in the above combinatorial
interpretation, one may
replace $y$ with $-y$.
For $n\geq 2$, we come to the determinant
 \begin{align*}
  V_n(x,y) = \begin{vmatrix}
  x &  1&  &  &  &  \\
  -y &  x &  2 &  & & \\
   & -(y+1) &  x &  \quad  \ddots & & \\
  &  & \ddots  &  \ddots &  \ddots &   \\
  & &  &  \ddots &  x &   n-1\\
  &  &  &  &   -(y+n-2)&   x
\end{vmatrix}_{n \times n}.
\end{align*}
The initial values of $V_n(x,y)$ for $n=0,1$ remain the same as those of $U_n(x,y)$.
The first few values of $V_n(x,y)$ are given below:
\begin{eqnarray*}
    V_0(x,y) & = & 1, \\[3pt]
    V_1(x,y) & = & x, \\[3pt]
    V_2(x,y) & = & x^2+y, \\[3pt]
    V_3(x,y) & = & x^3+3xy+2x, \\[3pt]
    V_4(x,y) & = & x^4 + 6x^2 y +8 x^2+ 3y^2 + 6y ,\\[3pt]
    V_5(x,y) & = & x^5 + 10 x^3 y+20 x^3+15 x y^{2}+50 x y+24 x.
     \end{eqnarray*}

Multiplying each odd row and each odd column of the
determinant $V_n(x,y)$ by $-1$, it can be recast as
\begin{align} \label{V-n}
   V_n(x,y) &=\begin{vmatrix}
  x &  -1&  &  &  &  \\
  y &  x &  -2&  & & \\
   & y+1 &  x & \ddots & & \\
  &  & \ddots  &  \ddots &  \ddots &   \\
  & &  & \ddots &  x &   -(n-1)\\
  &  &  &  &   y+n-2&   x
  \end{vmatrix}_{n \times n}.
\end{align}
Indeed, the above form of the Cayley continuants is
the basis of our wide band extension.

The purpose of this paper is to
make a connection between the wide band
Cayley continuants and
the joint distribution of the number of
$r$-regular cycles and the number of $r$-singular cycles. A cycle is called $r$-regular if its length is not divisible by $r$, and is called $r$-singular if its length is divisible by $r$.
We show that the wide band Cayley continuants can be interpreted
as the polynomials for the joint distribution of
the foregoing statistics.

\section{The wide band Cayley continuants }

We introduce the wide band Cayley continuants
as an extension of the original Cayley continuants
based on the expression of $V_n(x,y)$ in (\ref{V-n}).
Note that the term continuant is kept in a larger sense,
since the wide band Cayley continuants are not restricted to
be tridiagonal. For the sake of rigor,   the  wide
band Cayley determinants may be considered
as an alternative terminology.

Recall that an $n\times n$ matrix $A=(a_{i,\,j})_{n \times n}$ is
called a $(p,q)$-band matrix if $a_{i,\,j} = 0$ for
$j > i + p$ or $i > j + q$, where $0 \leq p,q \leq n-1$,
see B\"orgers \cite[Definition 5.1]{Borgers-2022}. Like the case of $r=2$, we use
$V_n^{(r)}(x,y)$ to denote the
wide band Cayley continuants, which generalize the Cayley continuants $U_n(x,y)$ with $y$ replaced by $-y$.

Consider a specific kind of $n \times n$ $(1,r-1)$ band matrices $A_n^{(r)}(x,y)$, whose determinants are called the wide band Cayley continuants, denoted by $V_n^{(r)}(x,y)$. For
$r\geq 2$ and $n \geq r$, the first superdiagonal
(above the main diagonal) is
$\left(-1,-2,\cdots,-(n-1)\right)$,
the main diagonal and the $(r-2)$
subdiagonals (below the main diagonal) have the same entry $x$,
and the $(r-1)$-st subdiagonal is $(y, y-1,y+2,\cdots,y-n+r)$,
as illustrated below:
\begin{align*}
V_n^{(r)}(x,y)=\begin{vmatrix}
  x&   -1&  &  &  &  & & & \\
  \vdots &  x&   -2&  &  &  & & \\
  x&  \vdots &  x&  \ddots &  &  & & \\
   y&  x&  \vdots &  \ddots & \ddots &  & & \\
  &  y-1&  x&  \ddots &  \ddots &  \ddots  & & \\
  &  &  \ddots &  \ddots &  \ddots &   \ddots\;\;& \;\,\ddots & \\
  &  &  & \ddots \;\; & x&  \cdots &   x&   -(n-1) \\
   &  &  &  & y-n+r&  x&   \cdots &  x
\end{vmatrix}_{n \times  n}.
\end{align*}

In particular, for $r=2$ and $n \geq 2$,
$V_n^{(2)}(x,y)$ takes the form of $V_n(x,y)$ as in
(\ref{V-n}). For $ r =3 $ and $n=6$, the wide band Cayley continuant  $V_6^{(3)}(x,y)$ is given by
\begin{align*}
V_6^{(3)}(x,y)=\begin{vmatrix}
  x&   -1&  &  &  &    \\
  x &  x&   -2&  &  \\
  y &  x &  x&  -3 &    \\
 &  y-1&  x&  x &  -4 &   \\
 & &  y-2&  x&  x &  -5  \\
 & &  &  y-3 & x & x
\end{vmatrix} .
\end{align*}

As to the initial values, When $n=0$, set $A_0^{(r)}(x,y)$ to be the empty matrix
with determinant $V_0^{(r)}(x,y)=1$. When $n=1$, set $A_0^{(r)}(x,y)=(x)$ and $V_1^{(r)}(x,y)=x$.
For $2 \leq n \leq r-1$, we define
$V_n^{(r)}(x,y)$ to be a truncated form of the
general case, that is,
\begin{align} \label{V-i-2}
V_n^{(r)}(x,y)=\begin{vmatrix}
  x& -1&    &    &  &   \\
  x&  x&  -2&    &  &   \\
  x&  x&   x&  \ddots &  &   \\
  \vdots & \vdots &  \vdots &  \ddots & \ddots  &  \\
  x&  x&   x&  \cdots &  x&  -(n-1)\\
  x&  x&   x&  \cdots &  x&  x \\
\end{vmatrix}_{n \times  n}.
\end{align}
For example, we have
\begin{align*}
V_2^{(4)}(x,y)=\begin{vmatrix}
  x&  -1\\
  x&  x
\end{vmatrix} \qquad \text{and} \qquad V_3^{(4)}(x,y)=\begin{vmatrix}
  x&  -1&  \\
  x&  x&  -2\\
  x&  x&  x
\end{vmatrix}.
\end{align*}

The initial values of $V_{n}^{(r)}(x,y)$ for $0\leq n \leq r-1$ are given
as follows, where $x^{(n)}$ stands for the rising factorial,
that is, $x^{(0)}=1$, and for $n\geq 1$, \[ x^{(n)}=x(x+1)
\cdots (x+n-1).\]

\begin{lemma} \label{initial}
For $r \geq 2$ and $0\leq n \leq r-1$,
we have
\[ V_{n}^{(r)}(x,y) = x^{(n)}. \]
\end{lemma}

\proof For $n=0,1$, by definition there is nothing to be said.
For $n\geq 2$, expanding
the determinant in (\ref{V-i-2})
along the last column, we see that
$V_n^{(r)}(x,y)$ admits the
same recurrence relation as $x^{(n)}$.
\qed

Analogous to the recurrence relation for the
Cayley continuants,  the following
relation holds for $V_n^{(r)}(x,y)$, where
$(x)_n$ stands for the lower factorial, that is,
$(x)_0=1$ and for $n\geq 1$,
\[ (x)_n = x (x-1) \cdots (x-n+1). \]

\begin{theorem}
\label{recurrence}
    For $r \ge 2$ and $n \ge r$, we have
    \begin{align}
        V_n^{(r)}(x,y)=x\sum_{i=1}^{r-1}(n-1)_{i-1}\;V_{n-i}^{(r)}(x,y) +(y+n-r)(n-1)_{r-1}\;V_{n-r}^{(r)} (x,y).
    \end{align}
    where $V_{n}^{(r)}(x,y)=x^{(n)}$ for $0\leq n \leq r-1$.
\end{theorem}

\proof Let $A_n^{(r)}(x,y)=(a_{i,\,j})_{n \times n}$ denote the matrix of the determinant
$V_n^{(r)}(x,y)$. For $r \ge 2$  and $n \ge r$, expanding the determinant
$V_n^{(r)}(x,y)$ along the last row, we get
$$V_n^{(r)}(x,y)=\sum_{i=1}^{r}a_{n,\,n-i+1}C_{n,\,n-i+1},$$
where $C_{n,\,n-i+1}$ is the cofactor of $a_{n,\,n-i+1}$. For $i=1$, we have $$C_{n,\,n}=V_{n-1}^{(r)}(x,y) = (n-1)_{i-1}\;V_{n-i}^{(r)}(x,y).$$

Next, for $i=r$, consider the case when $a_{n,1}\not= 0$. It can
happen only when $n=r$. In this case, we have
$$C_{n,1}=(n-1)_{n-1}=(n-1)_{i-1}\;V_{n-i}^{(r)}(x,y),$$
where we have made use of the fact that
$V_0^{(r)}(x,y)=1$.

If $a_{n,\,1}=0$, that is, $n>  r$, for $2\leq i \leq r$, we have
\begin{align*}
    C_{n,\,n-i+1}=(-1)^{i-1}
\left|\begin{array}{c:cccc}
A_{n-i}^{(r)}(x,y) & 
\\
\hdashline
\text{\large *} & \begin{array}{cccc}	
    -(n-i+1) & & \\	
    &  -(n-i+2) &  
    \\
    &  &  \ddots &  \\
    \multicolumn{2}{c}{\raisebox{1.3ex}[0pt]{\large *}}
    & & -(n-1)
    \end{array}
\end{array}\right|.
\end{align*}
Taking
the signs into account, we see that
\begin{align*}
    C_{n,\,n-i+1}=(n-1)_{i-1}\,\left|A_{n-i}^{(r)}(x,y)\right|=(n-1)_{i-1}\,V_{n-i}^{(r)}(x,y).
\end{align*}
Since $a_{n,\,n-i+1}=x$ for $1 \le i \le r-1$ and
$a_{n,\,n-r+1}=y+n-r$, we obtain the required recurrence relation
for $V_n^{(r)}(x,y)$. By Lemma \ref{initial}, the initial values of $V_{n}^{(r)}(x,y)$ for $0\leq n \leq r-1$ are
given by $x^{(n)}$, as expected. \qed

We now turn to the exponential generating functions for
the wide band Cayley continuants. Write
\begin{align*}
    V_r(x,y;t)=\sum_{n \ge 0}V_n^{(r)}(x,y)\frac{t^n}{n!}.
\end{align*}
Theorem \ref{recurrence} implies the differential equation
\begin{align*}
    (1-t^r)V'_r(x,y;t)=\left(t^{r-1}y+x(1-t^{r-1})(1-t)^{-1}\right)V_r(x,y;t),
\end{align*}
with the initial value $V_r(x,y;0)=1$.
This gives the following formula for $V_r(x,y;t)$. As will be seen,
the combinatorial interpretation in the
next section along with the known generating functions also
leads to the same expression for $V_r(x,y;t)$ without
resorting to a differential equation.

\begin{theorem}
    For $r \ge 2$, we have
\begin{align}\label{V-rxy}
    V_r(x,y;t)=\left(1-t^r\right)^{\frac{x-y}{r}}(1-t)^{-x}.
\end{align}
\end{theorem}

For $r=2$,  the above formula reduces to
the generating function of $V_n(x,y)$ as
derived by  Munarini and Torri \cite{MT-2005}, that is,
\begin{align}
    V_2(x,y;t)=(1-t^2)^{\frac{x-y}{2}}(1-t)^{-x}=\frac{(1+t)^{(x-y)/2}}{(1-t)^{(x+y)/2}}.
\end{align}

\section{A combinatorial interpretation}

In this section, we give a combinatorial interpretation of the above wide band Cayley continuants.

For $n \ge 0$, define
$W_n^{(r)}(x,y)$ to be the polynomial for the
joint distribution of the number of
$r$-regular cycles and the number of $r$-singular
cycles over permutations of $[n]$ with the assumption that $W_0^{(r)}(x,y)=1$. To be more specific, define
\begin{align*}
    W_n^{(r)}(x,y)=\sum_{\sigma \in S_n}x^{r(\sigma)}y^{s(\sigma)},
\end{align*}
where  $r(\sigma)$ and $s(\sigma)$ denote
the number of $r$-regular cycles and the number of
$r$-singular cycles of $\sigma$, respectively.
For $r\geq 2$, by definition, we have
$W_1^{(r)}(x,y)=x$.

For $r=2$, Munarini and  Torri \cite{MT-2005} found
a combinatorial interpretation of $V_n(x,y)$, that is,
$$V_n(x,y)=W_n^{(2)}(x,y).$$

In the following theorem,
we give a recurrence relation for $W_n^{(r)}(x,y)$,
which turns out to coincide with that for the
wide band Cayley continuants.  The argument is along
the line of Munarini and  Torri \cite{MT-2005}.

\begin{theorem}
For $r \ge 2$ and $n \ge r$, we have
\begin{align}
    \label{recurrence-weight}
    W_{n}^{(r)}(x,y)=x\sum_{i=1}^{r-1}(n-1)_{i-1}\;W_{n-i}^{(r)}(x,y) +(y+n-r)(n-1)_{r-1}\;W_{n-r}^{(r)} (x,y),
\end{align}
where
$W_{n}^{(r)}(x,y)= x^{(n)}$ for $0 \leq n \leq r-1$.
\end{theorem}

\proof Let $\sigma$ be a permutation of $[n]$.
Assume $\sigma$ is represented in the cycle notation,
and each cycle begins with its minimum element.
In particular, we call the cycle containing the element
$1$ the first cycle.
There are three cases depending on the length of the first
cycle.

\noindent
Case 1: The first cycle length $i$ is less than $r$.
Then it is a $r$-regular cycle with weight $x$. Moreover,
there are  $(n-1)_{i-1}$ choices for the first cycle.
Define the weight of $\sigma$ to be
$x^{r(\sigma)}y^{s(\sigma)}$. Then
the sum of weights of all possible permutations  in this case
equals
$$x \sum_{i=1}^{r-1} (n-1)_{i-1}\, W_{n-i}^{(r)}(x,y).$$

\noindent
Case 2: The first cycle length $i$ equals $r$.
Then it is a $r$-singular cycle of weight $y$.
Moreover, there are $(n-1)_{r-1}$ choices for
the first cycle. The sum
of weights of all possible permutations in this case equals
$$y (n-1)_{r-1}\, W_{n-r}^{(r)}(x,y). $$

\noindent
Case 3:
The first cycle length $i$ is greater than  $r$.
Let $C=(1~j_2~\cdots~j_r~j_{r+1}~\cdots)$ be
the first cycle, and let $D=j_{r+1}~j_{r+2}~\cdots$ be
the permutation obtained from $C$ by removing the
first $r$ elements. Therefore, $\sigma$ can be
recovered from the cycle $C'=(1~j_2~\cdots~j_r)$,
the permutation $D$, and the rest of the cycles of
$\sigma$. Now, there are $(n-1)_{r-1}$ choices
for $C'$.

Assume that $C'$ is given.
There are $n-r$ elements left.
For any permutation $\pi$ on these remaining elements,
we need to specify an element
to play the role of $j_{r+1}$. There are
$n-r$ choices. Hence there are
$(n-1)_r$ choices for $j_2~j_3~\cdots~j_{r+1}$.
At this point, we may reconstruct
$C$ from $C'$, the specified element and the permutation $\pi$. Observe that
$C$ is $r$-regular if and
only if the cycle $(j_{r+1}\cdots)$ in $\pi$ is $r$-regular, and so
$\sigma$ and $\pi$ have the same weight.
Thus the sum of weights of all possible
permutations $\sigma$ in this case equals
$$(n-1)_{r} \, W_{n-r}^{(r)}(x,y).$$

Finally, when $1 \le n \le r-1$, each cycle in the permutation of $[n]$ is $r$-regular and has weight $x$.
By the combinatorial interpretation of the signless
Stirling numbers of the fist kind, we see that
$W_n^{(r)}(x,y)$ is precisely $x^{(n)}$. This completes the
proof.
\qed

In closing, we remark that in view of the
combinatorial interpretation of $V_n^{(r)}(x,y)$,
their generating function
can be deduced from the
generating functions for
$r$-regular permutations and
$r$-cycle permutations. A permutation is called
  $r$-regular if all its  cycles are $r$-regular, whereas
an $r$-cycle permutation is referred to a permutation
in which every cycle is $r$-singular.
For $n\geq 1$, let
${\rm Reg}_r(n)$ and ${\rm Cyc}_r(n)$ denote the set of all $r$-regular permutations and the set of all $r$-cycle permutations in $S_n$.
As usual, we set $|{\rm Reg}_r(0)|=1$ and $|{\rm Cyc}_r(0)|=1$.
Notice that  the notations $NODIV_r(n)$ and $PERM_r(n)$, are used
in  B\'ona-Mclennan-White \cite{BMW-2000}
in lieu of ${\rm Reg}_r(n)$ and ${\rm Cyc}_r(n)$.

  The following generating functions have
  long been known, see \cite{Bolker-Gleason-1980, BMW-2000, Erdos-Turan-1967}.
\begin{align*}
    \sum_{n \ge 0}|{\rm Reg}_r(n)|\frac{t^n}{n!} &=\exp\left(\sum_{n \neq 0 \bmod r}(n-1)!\frac{t^n}{n!}\right)
    =\frac{(1-t^r)^{1/r}}{1-t},\\[9pt]
    \sum_{n \ge 0}|{\rm Cyc}_r(n)|\frac{t^n}{n!} &=\exp\left({\sum_{\substack{n \ge 1, n = 0 \bmod r}}}(n-1)!\frac{t^n}{n!}\right)=(1-t^r)^{-1/r}.
\end{align*}
It follows that
 \begin{align*}
      V_r(x,y;t) &=
      \exp\left(x \sum_{n \neq 0 \bmod r}(n-1)!\frac{t^n}{n!}\right)
     \, \exp\left({y \sum_{\substack{n \ge 1, n = 0 \bmod r}}}(n-1)!\frac{t^n}{n!}\right) \\[9pt]
     & =
      \left(\sum_{n \ge 0}|{\rm Reg}_r(n)|\frac{t^n}{n!}\right)^x\left(\sum_{n \ge 0}|{\rm Cyc}_r(n)|\frac{t^n}{n!}\right)^y\\[9pt]
      &= \left(1-t^r\right)^{\frac{x-y}{r}}(1-t)^{-x},
 \end{align*}
 which is in agreement with (\ref{V-rxy}).

\vskip 6mm
\noindent{\large\bf Acknowledgment.}
 This work was supported
by the National
Science Foundation of China.


\begin{thebibliography}{99}

\bibitem{Bolker-Gleason-1980}
E.D. Bolker and A.M. Gleason,
Counting permutations,
J. Combin. Theory Ser. A, 29 (1980) 236--242.

\bibitem{BMW-2000}
M. B\'ona, A. Mclennan, and D. White,
Permutations with roots,
Random Struct. Algor., 17 (2000) 157--167.

\bibitem{Borgers-2022}
C. B\"orgers,
Introduction to Numerical Linear Algebra,
SIAM, 2022.

\bibitem{Cayley-1858}
A. Cayley,
On the determination of the value of a certain determinant,
Quart. J. Math., {\romannumeral2} (1858) 163--166 (Collected Mathematics Papers, Vol. 3, Cambridge University press, Cambridge, 1919, pp. 120--123).

\bibitem{Erdos-Turan-1967}
P. Erd\H{o}s and P. Tur\'an,
On some problems of a statistical group-theory. {\uppercase\expandafter{\romannumeral2}},
Acta Math. Acad. Sci. Hungar., 18 (1967) 151--163.

\bibitem{Munarini-2022}
E. Munarini,
Umbral operators for Cayley and Sylvester continuants,
Appl. Anal. Discrete Math., 16 (2022) 307--327.

\bibitem{MT-2005}
E. Munarini and D. Torri,
Cayley continuants,
Theoret. Comput. Sci., 347 (2005) 353--369.


\end{thebibliography}
\end{document}